\documentclass[11pt]{article}
\usepackage{amsfonts}

\usepackage{amsfonts, amsmath, amssymb}
\usepackage{amssymb,amsfonts,amsmath,
latexsym, epsfig,cite, psfrag,eepic,colordvi}
\usepackage{amscd,graphics}
\usepackage{lineno}

\newcommand{\qed}{\hfill $\Box $}
\newcommand{\pf}{\noindent {\bf Proof.} }

\newtheorem{theorem}{Theorem}[section]
\newtheorem{lemma}{Lemma}[section]
\newtheorem{coro}[theorem]{Corollary}

\newtheorem{conjecture}[theorem]{Conjecture}

\begin{document}

\title{$L$-factors and  adjacent vertex-distinguishing edge-weighting}

\date{}

\author{Yinghua Duan\,\textsuperscript{a}, Hongliang Lu\,\textsuperscript{b}$\thanks{  {\it Corresponding
email (Lu): luhongliang215@sina.com}}$,  \ Qinglin
Yu\,\textsuperscript{c} \
\\ {\small \textsuperscript{a}Center for Combinatorics, LPMC}
\\ {\small Nankai University, Tianjin, China}
\\ {\small \textsuperscript{b} Department of Mathematics}
\\ {\small Xi'an Jiaotong University, Xi'an,  China}
\\ {\small \textsuperscript{c}Department of Mathematics and Statistics}
\\ {\small Thompson Rivers University, Kamloops, BC, Canada}
}

\maketitle

\begin{abstract}
An edge weighting problem of a graph $G$ is an assignment of an
integer weight to each edge $e$. Based on edge weighting problem,
several types of vertex-coloring problems are put forward. A simple
observation illuminates that edge weighting problem has a close
relationship with special factors of graphs. In this paper, we
obtain several results on the existence of factors with the
pre-specified degrees, which generalizes earlier results in
\cite{ADMRT, ADR}. Using these results, we investigate
edge-weighting problem. In particular, we prove that every
$4$-colorable graph admits a vertex-coloring $4$-edge-weighting.
\end{abstract}

\section{Introduction}
In this paper, we consider only finite, undirected and simple
graphs. For a graph $G=(V,E)$, if $v\in V(G)$ and $e\in E(G)$, we
use $v\sim e$ to denote that $v$ is an end-vertex of $e$. For $v\in
V(G)$, $N_G(v)$ denotes the set of vertices which are adjacent to
$v$. For a spanning subgraph $H$ of $G$ and $W \subseteq V(G)$, we
use $d_H(v)$ for the number of neighbors of $v$ in $H$ and
$d_W(v)=|N_G(v)\cap W|$. In addition, we use $\omega(H)$ to denote
the number of connected components of $H$. A $k$-{\em vertex
coloring} $c$ of $G$ is an assignment of $k$ integers, $1,2,\dots,
k$, to the vertices of $G$, and the color of a vertex $v$ is denoted
by $c(v)$. The coloring is {\em proper} if no two adjacent vertices
share the same color. A graph $G$ is $k$-{\em colorable} if $G$ has
a proper $k$-vertex coloring. The {\em chromatic number} $\chi(G)$
is the minimum number $r$ such that $G$ is $r$-colorable. Notations
and terminologies that are not defined here may be found in
\cite{BB}.

A $k$-{\it edge-weighting}  of a graph $G$ is an assignment $w:\
E(G)\rightarrow \{1,\dots , k\}$. An edge weighting naturally
induces a vertex coloring $c(u)$ by defining $c(u)=\sum_{u\sim e}
w(e)$ for every vertex $u \in V(G)$. A $k$-edge-weighting of a graph
$G$ is \textit{vertex-coloring} if the induced vertex-coloring is
proper, that is, $c(u)\neq c(v)$ and say that $G$ admits a
\emph{vertex-coloring $k$-edge-weighting}.

A $k$-edge-weighting  can also be viewed as a partition of edges
into sets $\{S_1, \cdots , S_k\}$. For each vertex $v$, let $X_v$
denote the multiset in which the elements are the weightings of the
edges incident with $v$ and the multiplicity of $i$ in $X_v$ is the
number of edges incident to $v$ in $S_i$. An edge-weighting is
\textit{proper} if no two incident edges receive the same label. An
edge-weighting is \textit{adjacent vertex-distinguishing} if for
every edge $e=uv$, $X_u\neq X_v$; it is \emph{vertex-distinguishing}
if $X_u\neq X_v$ holds for any pair of vertices $u,v\in V(G)$.
Proper (adjacent) vertex-distinguishing edge-weighting has been
studied by many researchers (see \cite{ATT, BRS,BS}) and is
reminiscent of harmonious colorings (see \cite{Edwards}). Clearly,
if a $k$-edge-weighting is vertex-coloring, then it is adjacent
vertex-distinguishing. However, the converse may not hold.

If a graph have an edge as a component, it cannot have an adjacent
vertex-distinguishing or vertex-coloring edge-weighting. So in this
paper, we only consider the graphs without an edge component, we
refer such graphs as \textit{nice graphs}.

In \cite{KLT}, Karo\'{n}ski, {\L}uczak and Thomason initiated the
study of vertex-coloring and adjacent vertex-distinguishing
edge-weightings and  they brought forward the following conjecture.

\begin{conjecture}\rm{(Karo\'{n}ski, {\L}uczak and Thomason, \cite{KLT})}
Every nice graph admits a vertex-coloring 3-edge-weighting.
\end{conjecture}

Furthermore, they proved that the conjecture holds for 3-colorable
graphs (see Theorem 1 of \cite{KLT}). In \cite{CLWY}, Chang
\textit{et al.} considered bipartite graphs $G=(X,Y)$ and proved
that if  $|X||Y|$ is even, then the graph admits a vertex-coloring
$2$-edge-weighting. In \cite{DLYZ}, Duan \textit{et al.} improve
this result and show that all 3-connected bipartite graphs have
vertex-coloring 2-edge-weighting. For general graphs, Addario-Berry
\textit{et al.} showed that every nice graph admits a
vertex-coloring 30-edge-weighting \cite{ADMRT}. In \cite{ADR},
Addario-Berry, Dalal and Reed improved this result and proved that
every nice graph permits a vertex-coloring 16-edge-weighting.
 Wang and Yu \cite{WY} improved this bound to 13. Recently, Kalkowski, Karo\'nski and
Pfender \cite{Kalkowski} proved every nice graph permits a
vertex-coloring 5-edge-weighting.

On the other hand, there are many results  for adjacent
vertex-distinguishing edge-weighting. In \cite{KLT}, it is showed
that every nice graph permits an adjacent vertex-distinguishing
213-edge-weighting and graphs with minimum degree at least $10^{99}$
permit an adjacent vertex-distinguishing 30-edge-weighting. In
\cite{AADR}, it is showed that every nice graph permits an adjacent
vertex-distinguishing 4-edge-weighting and that graphs of minimum
degree at least 1000 permit an adjacent vertex-distinguishing
3-edge-weighting.

For a graph $G$, there is a close relationship between
$2$-edge-weighting and a special factor. Let $L: V (G)\mapsto 2^N$
be a set function, a \emph{list factor} (or \emph{$L$-factor} for
short) of a graph $G$ is a spanning subgraph $H$ such that
$d_H(v)\in L(v)$ for all $v\in V(G)$.

In general, $L$-factor problem is a NP-complete problem, even when
$G$ is bipartite. The comprehensive investigation of $L$-factors was
carried out by Lov\'asz \cite{lov2}. So an $L$-factor is a spanning
subgraph with degrees from specified sets. When each $L(v)$ is an
interval, then $L$-factors are the same as the usual degree factors.
For instance,
 let $f$ and $g$ be nonnegative
integer-valued functions on $V(G)$ with $f\geq g$ and
$L(v)=[g(v),f(v)]$ for $v\in V(G)$, then an $L$-factor is exactly a
$(g,f)$-factor. It was  shown that every graph has a spanning
subgraph in which every vertex has pre-specified degree
\cite{ADMRT,ADR}. In Section 2, we generalize earlier results in
\cite{ADMRT, ADR} about $L$-factors. Using these results, in Section
3, we show that nice graphs with appropriate degree condition have
an adjacent vertex-distinguishing 2-edge-weighting.

Weighting the edges of a graph with elements of a group $\Gamma$
also gives rise to a vertex-coloring. If the vertex-coloring is
proper, we say that $G$ admits a \textit{vertex-coloring
$\Gamma$-edge-weighting}. The edge-weighting problem on groups has
been studied by Karo\'{n}ski, {\L}uczak and Thomason in \cite{KLT}.
They proved that if $\Gamma$ is a finite abelian group of odd order
and $G$ is a non-trivial $|\Gamma|$-colorable graph, then $G$ admits
a vertex-coloring $\Gamma$-edge-weighting. In Section 4, we obtain
several results on vertex-coloring $\Gamma$-edge-weighting. Using
these results, we deduce that every $4$-colorable graph admits a
vertex-coloring $4$-edge-weighting.

\section{Subgraphs with pre-specified degree}

The following results about $L$-factors were proved in \cite{ADMRT}
and \cite{ADR}, respectively. In this section, we generalize these
results.

\begin{theorem} \rm{(Addario-Berry {\it et al.}, \cite{ADMRT})}\label{factor2}
Let $G$ be a graph and $L(v)=\{a_{v}^{-}, a_{v}^{-}+1, a_{v}^{+},
a_{v}^{+}+1\}$, for every $v\in V(G)$  such that
$\frac{d_{G}(v)}{3}\leq a_{v}^{-} \leq\frac{d_{G}(v)}{2}\leq
a_{v}^{+}\leq\frac{2d_{G}(v)}{3}$. Then $G$ contains an $L$-factor.
\end{theorem}

\begin{theorem}\rm{(Addario-Berry, Dalal and Reed, \cite{ADR})} \label{Berry}
Let $G$ be a graph and $L(v)=\{a_{v}^{-}, a_{v}^{-}+1, a_{v}^{+},
a_{v}^{+}+1\}$, for every $v\in V(G)$  such that $a_{v}^{-}\leq
\lfloor \frac{1}{2}d(v)\rfloor\leq a_{v}^{+}<d(v)$, and
\begin{align*}
a_{v}^{+}\leq \min\{\frac{1}{2}(d(v)+a_{v}^{-})+1,
2(a_{v}^{-}+1)+1\}.
\end{align*}
 Then $G$ contains an $L$-factor.
\end{theorem}

\begin{theorem} \rm{(Addario-Berry, Dalal and Reed, \cite{ADR})} \label{lem7}
Let  $G=(X,Y)$ be bipartite graph. For $v\in X$, let $a_v^-=\lfloor
d(v)/2 \rfloor$ and  $a_v^+=a_v^-+1$. For $v\in Y$, choose $a_v^-$,
$a_v^+$ such that $a_v^-\leq d(v)/2 \leq a_v^+$ and $a_v^+\leq
\min\{\frac{d(v)+a_v^-}{2}+1,2a_v^-+1\}$. Let $L(v)=\{a_{v}^{-},
a_{v}^{+}\}$ for every vertex $v\in V$. Then $G$ contains an
$L$-factor.
\end{theorem}

The following characterization of $(g,f)$-factors was given by
Heinrich \emph{et al.} \cite{H} for the case of $g<f$.

\begin{theorem}\rm{(Heinrich {\it et al.}, \cite{H})} \label{Heinrich}
Let $G=(V,E)$ be a graph and, for all $v\in V$, integers $a_v$,
$b_v$ such that $0\leq a_v\leq b_v\leq d(v)$. Assume that one of the
following two conditions holds:
\begin{enumerate}
\item[(a)] $a_{v} < b_{v}$ for all $x \in V$;
\item[(b)] $G$ is bipartite.
\end{enumerate}
Then $G$ has a ($a_{v}, b_{v}$)-factor if and only if for all
disjoint sets of vertices of $A$ and $B$,

\begin{equation}\label{Heinrich-1}
\sum_{v\in A}(a_v-d_{G-B}(v))\leq \sum_{v\in B}b_v.
\end{equation}
\end{theorem}

\begin{theorem}\label{factor}
Suppose that $c_{1}$, $c_{2}$ and $c_{3}$ are three constants, where
$0<c_{1}<c_{2}<c_{3}<1$, $c_{3}-\frac{c_{1}(1-c_{3})}{1-c_{1}}\leq
c_{2}$ and $c_{1}\geq c_{3}c_{2}$. Let $G$ be a graph and
$L(v)=\{a_{v}^{-}, a_{v}^{-}+1, a_{v}^{+}, a_{v}^{+}+1\}$ for every
vertex $v\in V(G)$ such that $c_{1}d_{G}(v)\leq a_{v}^{-}\leq
c_{2}d_{G}(v)\leq a_{v}^{+}\leq c_{3}d_{G}(v)$. Then $G$ contains an
$L$-factor.
\end{theorem}

\pf Given a set of integers $S=\{a_{v}\ |\ v\in V\}$ and a subgraph
$H$ of $G$, we define the \emph{deficiency} of $H$ with respect to
 $S$ to be
\begin{align*}
def(H)=\sum_{v\in V(G)}\max\{0,a_{v}-d_{H}(v)\}.
\end{align*}
Suppose that $G$ contains no  $L$-factor. Choose $a_{v}\in
\{a^{-}_{v}, a^{+}_{v} \}, b_{v} = a_{v} + 1$ and a spanning
subgraph $H$ of $G$ such that for all $v\in V, d_{H}(v)\leq b_{v}$,
so that the deficiency is minimized over all such choices.
Necessarily, there exists at least one vertex $v\in V$ such that
$d_{H}(v)<a_{v}$, so the deficiency of $H$ is positive.

Let $A_{0}=\{v \ | \ d_{H}< a_{v}\}$. An $H$-alternating trail is a
trail $P=v_{0}v_{1}\ldots v_{k}$ with $v_{0}\in A_{0}$ and
$v_{i}v_{i+1}\not\in H$ for $i$ even, $v_{i}v_{i+1}\in H$ for $i$
odd. Let $$A=\{v \ | \  \mbox{there is an even $H$-alternating trail
ending at} \ v\},$$ and $$B=\{v \ | \   \mbox{there is an odd
$H$-alternating trail ending at}\ v\}.$$ Note that $A_{0}\subseteq
A$. For $v\in A$, $d_{H}(v)\leq a_{v}$, or else by alternating the
edges in $H$ along an even alternating trail ending in $v$, we
obtain a subgraph with less deficiency. Similarly, for $v\in B$,
$d_{H}(v)=b_{v}$ or else we can likewise decrease the deficiency by
alternating the edges  of $H$, this time along an odd alternating
trail ending at $v$. Note that $b_{v}>a_{v}$  implies that $A$ and
$B$ are disjoint. Furthermore, for any edge $e$  with one end in $A$
and other end not in $B$, then  $e\in E(H)$;  for any edge $e$ with
one end in $B$ and other end not in $A$, then $e\not\in E(H)$. By
these observations, we have
\begin{align}\label{main}
\sum_{v\in A}a_{v}&>\sum_{v\in A}d_{H}(v) =\sum_{v\in
B}d_{H}(v)+\sum_{v\in A}d_{G-B}(v)=\sum_{v\in
A}d_{G-B}(v)+\sum_{v\in B}b_{v},
\end{align}
which implies that (a) of Theorem \ref{Heinrich} fails for the sets
$A$ and $B$.

We make two claims:
\begin{align}\label{claim11}
\mbox{For all} \ v\in A, a_{v}-d_{G-B}(v)\leq c_{2}d_{B}(v)
\end{align}
and
\begin{align}\label{claim21}
\mbox{for all} \ v\in B, b_{v}\geq c_{2}d_{A}(v).
\end{align}
Then (\ref{claim11}) and (\ref{claim21}) together with the fact
$\sum_{v\in A}d_{B}(v)=\sum_{v\in B}d_{A}(v)$ imply inequality
(\ref{Heinrich-1}) holds for the sets $A$ and $B$,  a contradiction
to (\ref{main}). So to end the proof, we only need to prove
(\ref{claim11}) and (\ref{claim21}).

To see (\ref{claim11}), we consider $v\in A$. Assume that
$d_{H}(v)<a_{v}$. (Note that alternating  edges in $H$ along an even
alternating trail does not change the deficiency and the sets $A$
and $B$. Thus we  assume that any  vertex $v\in A$ satisfies
$d_{H}(v)<a_{v}$.) Furthermore, we may assume $a_{v} =a_{v}^{+}\geq
c_{2}d(v)$ or else (\ref{claim11}) holds automatically. We may also
assume that $d_{G-B}(v)>a_{v}^{-}+1$, otherwise, by setting
$a_{v}=a_{v}^{-}$ and removing  some edges  in $H$ from $v$ to $B$,
we can reduce the deficiency. Thus,
\begin{align*}
d_{G-B}(v)=d_{H-B}(v)>a_{v}^{-}+1\geq c_{1}d(v)
\end{align*}
which implies
$d_{B}(v)<(1-c_{1})d(v),$
and hence
$d_{G-B}(v)>\frac{c_{1}}{(1-c_{1})}d_{B}(v).$
So,
\begin{align*}
a_{v}-d_{G-B}(v)&\leq c_{3}d(v)-d_{G-B}(v)\\
&=c_{3}d_{B}(v)-(1-c_{3})d_{G-B}(v)\\
&<c_{3}d_{B}(v)-(1-c_{3})\frac{c_{1}}{(1-c_{1})}d_{B}(v)\\
&\leq c_{2}d_{B}(v).
\end{align*}

Next we show (\ref{claim21}). Let $v\in B$. We may assume that
$b_{v}=a_{v}^{-}+1<c_{2}d(v)$ or else (\ref{claim21}) holds
trivially. Suppose  that the statement fails, then
$d_{A}(v)>b_{v}/c_{2}\geq (c_{1}d(v)+1)/c_{2}\geq
(c_{1}d(v))/c_{2}+1/c_{2}> a_{v}^{+}$. There are $d_{A}(v)-b_{v}$
edges from $v$ to $A$ that are not in $H$. In particular, there is a
$w\in N(v)\cap A$ such that $vw\not\in H$. As noted above, we can
ensure that $d_{H}(w)<a_{w}$. This will not change the fact that
$vw\not\in H$. Setting $a_{v}=a^{+}_{v}$ and adding
$a_{v}^{+}-d_{H}(v)$ edges from $v$ to $A$ into $H$ (including the
edge $vw$), we decrease the deficiency. \qed

\vspace{3mm}

\textbf{Remark:} Theorem \ref{factor} is a generalization of Theorem
\ref{factor2}. To see this, let $c_{1}=\frac{1}{3},
c_{2}=\frac{1}{2},c_{3}=\frac{2}{3}$, then $c_{1},c_{2}$ and
 $c_{3}$ satisfy the conditions in Theorem \ref{factor}.

The following result is another extension of Theorem \ref{factor2}.
This time, we consider three consecutive pairs in $L(v)$.

\begin{theorem}\label{factor3}
For every vertex $v$ of graph $G$, suppose that we have chosen three
integers $a_{v}^{1},a_{v}^{2},a_{v}^{3}$ such that
$\frac{3}{10}d_{G}(v)\leq a_{v}^{1} \leq\frac{4}{10}d_{G}(v)\leq
a_{v}^{2}\leq\frac{6}{10}d_{G}(v)\leq a_{v}^{3}\leq
\frac{7}{10}d_{G}(v)$. Let $L(v)=\{a_{v}^{1}, a_{v}^{1}+1,
a_{v}^{2}, a_{v}^{2}+1,a_{v}^{3},a_{v}^{3}+1\}$  for every $v\in
V(G)$. Then $G$ contains an $L$-factor.
\end{theorem}

\pf Suppose the theorem doesn't hold. By Theorem \ref{factor},
choose $a_{v}\in \{a^{1}_{v}, a^{2}_{v}, a^{3}_{v}\}, b_{v} = a_{v}
+ 1$ and a spanning subgraph $H$ of $G$ with $d_{H}(v)\leq b_{v}$
for all $v\in V$, so that the deficiency is minimized over all such
choices. 
We construct $A$ and $B$ as in the proof of Theorem \ref{factor}. We
shall prove the following two claims:
\begin{align}\label{claim26}
\mbox{for all} \ v\in A, a_{v}-d_{G-B}(v)\leq d_{B}(v)/2
\end{align}
and
\begin{align}\label{claim2}
\mbox{for all} \ v\in B, b_{v}\geq d_{A}(v)/2
\end{align}

If (\ref{claim26}) and (\ref{claim2}), then we have
\begin{align*}
\sum_{v\in A}(a_{v}-d_{G-B}(v))\leq \frac{1}{2}\sum_{v\in
A}d_{B}(v)=\frac{1}{2}\sum_{v\in B}d_{A}(v)\leq \sum_{v\in B}b_{v},
\end{align*}
a contradiction. So it remains to prove (\ref{claim26}) and
(\ref{claim2}).

To see (\ref{claim26}), consider $v\in A$, and assume
$d_{H}(v)<a_{v}$. We may assume $a_{v}\in \{a_{v}^{2},a_{v}^{3}\}$
or else (\ref{claim26}) holds trivially. If $a_{v}=a_{v}^{3}$, then
we assume that $d_{G-B}(v)>a_{v}^{2}+1$, or else by letting
$a_{v}=a_{v}^{2}$ and removing from $H$ some of the edges from $v$
to $B$, we can reduce the deficiency. Moreover,
$d_{H-B}=d_{G-B}>a_{v}^{2}+1$, as, otherwise letting
$a_{v}=a_{v}^{2}$ and deleting edges of $H$ between $v$ and $B$
contradicts our choice of $H$. Thus
\begin{align*}
d_{H-B}=d_{G-B}>a_{v}^{2}+1\geq \frac{4}{10}d_{G}(v),
\end{align*}
which implies $d_{B}(v)<\frac{3}{5}d(v)$ and hence
$\frac{2}{3}d_{B}(v)<d_{G-B}(v).$ So we have
\begin{align*}
a_{v}-d_{G-B}(v)&\leq \frac{7}{10}d(v)-d_{G-B}(v)\\
&=\frac{7}{10}d_{B}(v)-\frac{3}{10}d_{G-B}(v)\\
&<\frac{7}{10}d_{B}(v)-\frac{3}{10}\ast\frac{2}{3}d_{B}(v)\\
&\leq \frac{1}{2}d_{B}(v).
\end{align*}
With similar arguments, (\ref{claim26}) holds for the case of
$a_{v}=a_{v}^{2}\geq \frac{1}{2}d_{G}(v)$.

Next, we show (\ref{claim2}). Let $v\in B$. We may assume that
$b_{v}=a_{v}^{i}+1$, where $i=1\ \mbox{or}\ 2$, otherwise
(\ref{claim2}) holds trivially. Suppose $b_{v}=a_{v}^{1}+1$ and
$d_{A}(v)>2b_{v}$, then there are $d_{A}(v)-b_{v}$ edges from $v$ to
$A$ that are not in $H$, in particular, there is a vertex $w\in
N(v)\cap A$, $vw\not\in H$. As noted above, we can ensure that
$d_{H}(w)<a_{w}$. This does not change the fact that $vw\not\in H$.
Setting \textbf{$a_{v}=a^{2}_{v}$} and adding $a_{v}^{2}-d_{H}(v)$
edges from $v$ to $A$ into $H$ (including the edge $vw$), we
decrease the deficiency. If \textbf{$a_{v}=a_{v}^{2}$}, the
arguments are similar. \qed

Using similar proof with a slight change, we obtain the following
two results.
\begin{theorem}\label{factor4}
Let $G$ be a graph and $c$ be a constant satisfing $0<c<2/3$. For
all $v\in V(G)$, given integers $a_{v}^{-}, a_{v}^{+}$ such that
$a_{v}^{-}\leq cd(v)\leq a_{v}^{+}<d(v)$, and
\begin{align}\label{inequality1}
a_{v}^{+}\leq \min\{cd(v)+(1-c)a_{v}^{-}+1, (a_{v}^{-}+1)/c+1\}.
\end{align}
Let $L(v)=\{a_{v}^{-},a_{v}^{-}+1,a_{v}^{+},a_{v}^{+}+1\}$. Then $G$
contains an $L$-factor.
\end{theorem}

\pf Similar to the proofs of Theorems  \ref{factor} and
\ref{factor3}, it is sufficient to prove the following two claims:
\begin{align}\label{claim3}
a_{v}-d_{G-B}(v)\leq cd_{B}(v) \ \ \mbox{for all} \ v\in A
\end{align}
and
\begin{align}\label{claim4}
b_{v}\geq cd_{A}(v) \ \ \mbox{for all} \ v\in B.
\end{align}
These two statements together with the fact that $\sum_{v\in
A}d_{B}(v)=\sum_{v\in B}d_{A}(v)$ imply (1) holds for the sets $A$
and $B$, completing the proof of Theorem \ref{factor4} by
contradiction.

To see (\ref{claim3}), consider $v\in A$ and assume
$d_{H}(v)<a_{v}$. We may assume $a_{v}=a_{v}^{+}$ or else
(\ref{claim3}) holds trivially. If $a_{v}=a_{v}^{-}$, then we assume
that $d_{G-B}(v)>a_{v}^{-}+1$, or else by letting $a_{v}=a_{v}^{-}$
and removing from $H$ some of the edges from $v$ to $B$, we can
reduce the deficiency. Moreover, $d_{H-B}=d_{G-B}>a_{v}^{-}+1$. Now
we have
\begin{align*}
a_{v}&\leq  cd(v)+(1-c)a_{v}^{-}+1\\
&\leq cd(v)+(1-c)(d_{G-B}(v)-2)+1\\
&= cd_{B}(v)+d_{G-B}(v)+2c-1.
\end{align*}
Since $a_{v}$ is an integer and $2c-1<1/3$, then (\ref{claim3})
holds.

To prove (\ref{claim4}), consider any $v\in B$. We may assume $a_{v}
= a_{v}^{-}<cd(v)$, otherwise, the statement holds trivially.
Suppose that the statement fails, then $d_{A}(v)>(a_{v} + 1)/c$,
thus $d_{A}(v)\geq a^{+}_{v}$ by (\ref{inequality1}). There are
$d_{A}(v)-b_v$ edges from $v$ to $A$ that are not in $H$. In
particular, there is a vertex $w\in N(v)\cap A, vw\not\in H$. As
noted above, we can ensure that $d_{H}(w)<a_{w}$, which does not
change the fact that $vw\not\in H$. Setting $a_{v} =a^{+}_{v}$ and
adding $a^{+}_{v}-d_{H} (v)$ edges from $v$ to $A$ into $H$
(including the edge $vw$), then $def(H)$ decreases. \qed

\vspace{3mm}

\begin{theorem}\label{thm7}
Let $c$ be a constant with $0<c<2/3$ and $G=(X,Y)$ be a bipartite
graph. For $v\in X$, let $a_v^-=\lfloor cd(v) \rfloor$ and
$a_v^+=a_v^-+1$. For $v\in Y$, choose $a_v^-$, $a_v^+$ such that
$a_v^-\leq \lfloor cd(v)\rfloor \leq a_v^+$ and $a_v^+\leq \min
\{cd(v)+(1-c)a_v^-+1,\frac{a_v^-}{c}+1\}$. Let
$L(v)=\{a_v^-,a_v^+\}$ for all $v\in V(G)$. Then $G$ contains an
$L$-factor.
\end{theorem}

\textbf{Remark:} It is easy to see that if we set $c=\frac{1}{2}$,
then Theorems \ref{factor4} and \ref{thm7}  becomes Theorems
\ref{Berry} and \ref{lem7}, respectively.


\section{Adjacent vertex-distinguishing edge-weighting}

It was proved  in \cite{KLT} that every $3$-colorable graph has a
vertex-coloring $3$-edge-weighting. In particular, it has an
adjacent vertex-distinguishing $3$-edge-weighting.  The natural
question is that if every $2$-colorable (bipartite) graph has an
adjacent vertex-distinguishing $2$-edge-weighting.

In \cite{CLWY}, Chang \emph{et al.} considered this problem and
obtained the following results.

\begin{lemma} \rm{(Chang {\it et al.}, \cite{CLWY})} \label{yu}
A non-trivial connected bipartite graph $G=(U,W)$ admits a
vertex-coloring $2$-edge-weighting if one of following conditions
holds:

    $(1)$ $|U|$ or $|W|$ is even;

    $(2)$ $\delta(G) =1$;

    $(3)$ $\lfloor d(u)/2 \rfloor +1 \neq d(v)$ for any edge $uv \in
    E(G)$.
\end{lemma}

Recently, Duan \emph{et al.} \cite{DLYZ} improved this result and
 proved that every 3-connected bipartite graph admits a
vertex-coloring 2-edge-weighting. Since a graph admits a
vertex-coloring $k$-edge-weighting implies that it has an adjacent
vertex-distinguishing $k$-edge-weighting. Next, we go on to study
adjacent vertex-distinguishing $2$-edge-weighting on bipartite
graph. We prove the following results.

\begin{theorem} Given a nice bipartite graph $G=(U,W)$. If there exists a vertex $v\in V(G)$ such that
$d_{G}(v)\not \in \{d_{G}(x)\ |\ x\in N(v)\}$, then $G$ admits an
adjacent vertex-distinguishing $2$-edge-weighting.
\end{theorem}

\pf If $|U|\cdot|W|$ is even, by Lemma \ref{yu}, the result is
clear. So we assume that both $|U|$ and $|W|$ are odd. Let $v\in U$
such that $d_{G}(v)\not \in \{d_{G}(x)\ |\ x\in N(v)\}$ . By Lemma
\ref{group2}, $G$ has a vertex-coloring $2$-edge-weighting, such
that $c(x)$  is odd for all $x\in U-v$ and  $c(y)$  is even for all
$y\in W\cup\{v\}$. Since $d_{G}(v)\not \in \{d_{G}(x)\ |\ x\in
N(v)\}$,  then $G$ admits an adjacent vertex-distinguishing
$2$-edge-weighting. Thus we complete the proof. \qed

\begin{theorem}\label{thm11} Every nice bipartite graph
with $\delta(G)\geq 6$ admits an adjacent vertex-distinguishing
2-edge-weighting.
\end{theorem}

\pf  Let $G=(U, W)$ be a bipartite graph. For $v\in U$, let
$a_v^-=\lfloor d(v)/2 \rfloor$ and $a_v^+=a_v^-+1$. For $v\in W$,
choose $a_v^-=\lfloor d(v)/2 \rfloor -1$ and $a_v^+=\lfloor d(v)/2
\rfloor +2$. Since $\delta(G)\geq 6$, in $W$, $a_v^-$ and $a_v^+$
satisfy the condition of Theorem \ref{lem7}, that is, $a_v^+\leq
min\{\frac{d(v)+a_v^-}{2}+1,2a_v^-+1\}$. So there is a spanning
subgraph $H$ such that $d_H(v)\in \{\lfloor d(v)/2 \rfloor, \lfloor
d(v)/2 \rfloor+1\}$ for all $v\in U$,  $d_H(v)\in \{\lfloor d(v)/2
\rfloor-1, \lfloor d(v)/2 \rfloor+2\}$ for $v\in W$. Thus, we can
label the edges in $E(H)$ with $1$ and edges in $G-E(H)$ with $2$,
which yield an adjacent vertex-distinguishing $2$-edge-weighting of
$G$. \qed

For non-bipartite graph,  Addario-Berry \emph{et al.} \cite{ADR}
proved that a graph $G$ with $\delta(G)\geq 12\chi(G)$ admits a
vertex-coloring $2$-edge-weighting. In \cite{DLYZ}, a lower bound of
minimum degree is improved to $8\chi(G)$, i.e., $\delta(G)\geq
8\chi(G)$, to ensure an adjacent vertex-distinguishing
$2$-edge-weighting in $G$.

\section{Vertex-coloring $\Gamma$-edge-weighting on graph}

In this section, we  consider the edge-weighting problem on groups.
The following technical lemmas are contained in \cite{KLT} and
\cite{DLYZ}, respectively.

\begin{lemma}\rm{(Karo\'{n}ski, {\L}uczak and Thomason, \cite{KLT})} \label{group1}
Let $\Gamma$ be a finite abelian group of odd order and  $G$  a
non-trivial $|\Gamma|$-colorable graph. Then there is a weighting of
the edges of $G$ with the elements of $\Gamma$ such that the induced
vertex weighting is proper coloring.
\end{lemma}

\begin{lemma}\rm{(Duan \emph{et al.}, \cite{DLYZ})} \label{group2}
Let $G$ be a connected nice  graph with chromatic number $k\geq 3$
and $\Gamma=\{g_1,g_2,\dots, g_k\}$ be a finite abelian group, where
$k=|\Gamma|$. Let $c_0$ be any $k$-vertex coloring of $G$ with color
classes $\{U_{1},\dots, U_k\}$, where $|U_i|=n_i$ for $1\leq i\leq
k$. If there exists an element $h\in \Gamma$ such that
$n_{1}g_{1}+\dots+n_{k}g_{k}=2h$, then there is an edge-weighting of
$G$ with the elements of $\Gamma$ such that the induced vertex
coloring is $c_0$.
\end{lemma}

Using Lemma \ref{group2}, we can prove the following result.

\begin{theorem}\label{mod4}
Let $Z_{r}$ with $r\equiv 0 \ (mod\ 4)$ be a cyclic group and $G$ be
a $r$-colorable graph.  Then there exists a vertex-coloring
$Z_{r}$-edge-weighting of  $G$.
\end{theorem}

\pf Let $\mathcal {U}: V(G)\mapsto Z_{r}$ be a proper color of $G$
with partition $(U_{1},\ldots, U_{r})$ and $\mathcal {U}(U_i)=i$. If
$\sum_{i=1}^{r}i|U_{i}|$ is even, then by Lemma \ref{group2}, the
result is followed. Now we assume $\sum_{i=1}^{r}i|U_{i}|$ is odd.
So $\sum_{i=1}^{\frac{r}{2}}(2i-1)|U_{2i-1}|$ is odd. Since $r/2$ is
even, so we can assume that there exists some $U_{2i-1}$ with even
order. If there exists some $U_{2j}$ with odd order, then we recolor
$U_{2j}$ with color $2i-1$ and $U_{2i-1}$ with $2j$ and the rest of
classes remains unchanged. Then $\sum_{i=1}^{r}\mathcal
{U}(i)|U_{i}|$ is even, by Lemma \ref{group2}, the result is
followed. Now  assume that $|U_{2k}|$ is even for $k=1,2,\ldots,
r/2$. Note that there exists a set $U_{2l-1}$ with odd order. Now we
recolor $U_{2k}$ with color $2l-1$ and $U_{2l-1}$ with $2k$ and the
rest of classes remains unchanged, then the result is followed. \qed

From the theorem above, the following results can be easily deduced.

\begin{theorem}
Let $G$ be a $4$-colorable graph. Then $G$ admits a vertex-coloring
$4$-edge-weighting.
\end{theorem}

\begin{coro}
Let $G$ be a $r$-colorable graph, where $r\neq 4k+2$. Then $G$ has a
vertex-coloring $r$-edge-weighting.
\end{coro}

Since every planar graph is 4-colorable, so we have the following
interesting result.

\begin{coro}
Every planar graph admits a vertex-coloring $4$-edge-weighting.
\end{coro}

\begin{theorem}\label{odd}
Let $G$ be a $r$-colorable graph. Suppose $G$ doesn't admit  a
vertex-coloring $Z_{r}$-edge-weighting. Let $\lambda: V(G)\mapsto
Z_{r}$ be an arbitrary proper color of $G$, then $|\lambda^{-1}(i)|$
is odd for $i=1,\ldots, r$.
\end{theorem}
\pf By Lemma \ref{group1} and Theorem \ref{mod4}, we can assume
$r\equiv 2$ (mod 4). Suppose the result doesn't hold. There exists a
set, say $\lambda^{-1}(i)$ with even order. Note that
$\sum_{l=1}^{r/2}(2l-1)|\lambda^{-1}(2l-1)|$ is odd. Thus there
exists some $l\neq i$ such that $|\lambda^{-1}(2l-1)|$ is odd. If
$i$ is even, then we recolor $\lambda^{-1}(2l-1)$ with   $i$ and
color $i$ with  $2l-1$ and obtain a coloring $\lambda'$. Then
$\sum_{l=1}^{r}l|\lambda'^{-1}(l)|$ is even, a contradiction to
Lemma \ref{group2}. So  $i$ is odd. Moreover, we can assume
$|\lambda^{-1}(2l)|$ is odd for $l=1,\ldots, r/2$. Now we recolor
$\lambda^{-1}(i)$ with  $2$ and $\lambda^{-1}(2)$ with  $i$, and
obtain a coloring $\lambda''$. Clearly
$\sum_{l=1}^{r}l|\lambda''^{-1}(l)|$ is even,  a contradiction
again. \qed

\begin{theorem}\label{lu} Let $G$ be a $k$-colorable graph, where
$(U_0,U_1,\ldots,U_{k-1})$ denote coloring classes of $G$. Then $G$
admits a vertex-coloring $k$-edge-weighting, if any of following
conditions holds:
\begin{itemize}
\item[$(i)$]$k\equiv 0 \pmod 4$;

\item[$(ii)$]$\delta(G)\leq k-2$;

\item[$(iii)$]there exists a class $U_{i}$ with $|U_{i}| \equiv
0 \pmod 2$ for some $i\in \{0,1,\ldots,k-1\}$;

\item[$(iv)$]$|V(G)|$ is odd.
\end{itemize}
\end{theorem}

\pf (i) By Lemma \ref{group1} and Theorem \ref{mod4}, the result is
followed.

(ii) Let $\lambda: V(G)\mapsto Z_{r}$ be a proper vertex coloring
with partition $(U_{1},\ldots,U_{r})$. By Theorem \ref{odd}, then
$|U_{i}|$ is odd for $i=1,\ldots,r$. Let $d_{G}(v)\leq r-2$ and
$v\in U_{i}$. Clearly, there exists some $U_{j}$ with $i\neq j$ such
that there are no edge between $v$ and $U_{j}$. We can recolor $v$
with  $j$ and the coloring of the rest vertices remain unchange.
Then we obtain a new coloring $\lambda'$, a contradiction to Theorem
\ref{odd}.

(iii) By Theorem \ref{odd}, the result is clear.

(iv) Consider $r\equiv 2$ (mod 4).   If $|G|=\sum_{i=1}^{r}|U_{i}|$
is odd, then there exists some $U_{i}$ such that $|U_{i}|$ is even.
We complete the proof.\qed

%
%
%
%

\end{document}